\documentclass[12pt]{article}
\usepackage[dvips]{graphicx}
\usepackage{amsmath,amsfonts,amssymb,amsthm,color,url}

\usepackage{natbib}

\newtheorem{theorem}{Theorem}

\newtheorem{corollary}{Corollary}

\newtheorem{proposition}{Proposition}

\newtheorem{definition}{Definition}
\newtheorem{example}{Example}

\newtheorem{problem}{Problem}

\usepackage{bm}
\bmdefine{\Bt}{t}
\bmdefine{\BX}{X}
\bmdefine{\BY}{Y}
\bmdefine{\BZ}{Z}
\bmdefine{\BB}{B}
\bmdefine{\BM}{M}
\bmdefine{\BD}{D}
\bmdefine{\Bi}{i}
\bmdefine{\Bj}{j}
\bmdefine{\Bx}{x}
\bmdefine{\By}{y}
\bmdefine{\Bz}{z}
\bmdefine{\Bv}{v}
\bmdefine{\Bw}{w}
\bmdefine{\Bn}{n}
\bmdefine{\Ba}{a}
\bmdefine{\Bb}{b}
\bmdefine{\Bc}{c}
\bmdefine{\Be}{e}
\bmdefine{\Bu}{u}
\bmdefine{\Bp}{p}
\bmdefine{\Bzero}{0}
\bmdefine{\Bone}{1}

\newcommand{\N}{{\mathbb N}}
\newcommand{\Z}{{\mathbb Z}}
\newcommand{\R}{{\mathbb R}}

\def\comment#1{\textit{[#1]}}
\def\comment#1{}

\title{Markov bases and subbases for bounded contingency tables}

\author{
Fabio Rapallo
\and
Ruriko Yoshida}
\date{}

\begin{document}

\maketitle

\begin{abstract}
In this paper we study the computation of Markov bases for
contingency tables whose cell entries have an upper bound. In
general a Markov basis for unbounded contingency table under a
certain model differs from a Markov basis for bounded tables.
\cite{Rapallo07} applied Lawrence lifting to compute a Markov
basis for contingency tables whose cell entries are bounded.
However, in the process, one has to compute the universal
Gr\"obner basis of the ideal associated with the design matrix for
a model which is, in general, larger than any reduced Gr\"obner
basis. Thus, this is also infeasible in small- and medium-sized
problems. In this paper we focus on bounded two-way contingency
tables under independence model and show that if these bounds on
cells are positive, i.e., they are not structural zeros, the set
of basic moves of all $2 \times 2$ minors connects all tables with
given margins.  We end this paper with an open problem that if we
know the given margins are positive, we want to find the necessary
and sufficient condition on the set of structural zeros so that
the set of basic moves of all $2 \times 2$ minors connects all
incomplete contingency tables with given margins.

{\bf keywords}: Structural zeros \and Markov basis \and Universal
Gr\"obner basis

\end{abstract}

\section{Introduction}

The study of statistical models to detect complex structures in
contingency tables has received great attention in the last
decades. (See \cite{agre:2002} for an overview of such models).
Among the main research themes in this field, here we consider
incomplete contingency tables (or equivalently, tables with
structural zeros) and models to go beyond independence in two-way
tables, such as quasi-independence models.

Contingency tables with upper bounds on the cell counts have
recently been considered in, e.g., \cite{cryan|dyer|randall:05}.
Bounded contingency tables can come, for instance, in the analysis
of designed experiments with multinomial response, as in
\cite{aoki|takemura:06}, and in logistic regression models, as in
e.g. \cite{chen-dinwoodie-dobra-huber2005}. We will use some
examples from these applications later in the paper.

In recent years, the use of algebraic and geometric techniques in
statistics has produced at least two relevant advances. One is a
better understanding of statistical models in terms of varieties
and polynomial equations, through the notion of toric models, as
described in Chapter 6 of \cite{pistone|riccomagno|wynn:01}.
Moreover, algebraic statistics has introduced a non-asymptotic
method for goodness-of-fit tests following a Markov Chain Monte
Carlo approach (see \cite{diaconis-sturmfels}). Such an algorithm
is based on the notion of Markov basis. In the last years the
computation of Markov bases for special statistical models has
involved both statisticians and algebraists.

In this paper we consider the computation of Markov bases for
bounded contingency tables. A general algorithm to compute Markov
bases for this case was described in \cite{Rapallo07}, using the
notions of Lawrence lifting and Universal Gr\"obner basis of a
polynomial ideal. When a Markov basis is computed through a
Universal Gr\"obner basis, we say that it is {\em Universal
Markov basis}. The Markov bases for these kind of tables are in
general very large, and we will show some explicit computations
later in the paper. Therefore the computation of smaller Markov
bases or subbases for special tables is a problem of major
interest.

In practice, computing the Markov basis for the bounded
contingency tables is infeasible because the number of elements in
the Markov basis is very large. However, for some cases, if we
know that the given margins are positive then the number of moves
connecting all tables is smaller than the number of elements in a
Markov basis for tables under the model.  Such connecting sets
were formalized in \cite{Chen06} with the terminology {\em Markov
subbases}. In this paper we consider bounded $I \times J$ tables
under independence model.  These tables are equivalent to $I
\times J \times 2$ tables under the models of no-3-way
interaction. Using this fact and the result from \cite{Ian2008},
in this paper, we show that if we know the bounds of cells are all
positive, that is, there are no structural zeros, then the set of
basic moves of all $2 \times 2$ minors connects all bounded
two-way contingency tables with given margins.

To summarize, we classify the  bounds of cells into the following
patterns:
\begin{enumerate}
\item[(i)] \label{case1} all cells are unbounded,

\item[(ii)] \label{case2} all cells are bounded by positive
integers,

\item[(iii)] \label{case3} some cells are unbounded and the others
are bounded by positive integers,

\item[(iv)] \label{case4} some cells are unbounded and the others
are structural zeros,

\item[(v)] \label{case5} some cells are bounded by positive
integers and the others are structural zeros,

\item[(vi)] \label{case6} all types of bounds appear.
\end{enumerate}
Case (i) is the standard case, already studied in
\cite{diaconis-sturmfels}. In the past,
\cite{aoki-takemura-2005jscs} dealt with the case (iv). In this
paper Theorem \ref{Univ_struct} deals with the case (v), Theorem
\ref{main} deals with  the case (ii), Section \ref{univ-part}
deals with the case (iii).

The organization of this paper is as follows. In Section
\ref{recalls} we recall the basic facts about Markov bases and
bounded contingency tables. In Section \ref{univ-inc} we present a
characterization of Universal Markov bases for incomplete tables,
showing that there is a simple connection between the Universal
Markov basis for an incomplete table and the corresponding
complete table. We present some explicit examples, focusing in
particular on quasi-independence models for two-way tables. In
Section \ref{univ-part} we show how to compute Markov bases when
the bounds involve only a subset of cell counts.  In Section
\ref{univ-subbases} we show our main theorem, that is, we consider
bounded two-way contingency tables under independence model. If we
know all bounds are positive (equivalently there are no structural
zeros), then the set of basic moves of all $2 \times 2$ minors
connects all bounded two-way contingency tables with given
margins. We end this paper with an open problem for incomplete
contingency tables with positive margins.

\section{Bounded contingency tables and Markov bases} \label{recalls}

Let ${\mathbf n}$ be a contingency table with $k$
cells. In order to simplify the notation, we denote by
$\mathcal{X}=\{1, \ldots, k \}$ the sample space of the
contingency table. In the special case of two-way tables with $I$
rows and $J$ columns, we will also denote the sample space with
${\mathcal X}=\{1, \ldots , I \} \times \{1 , \ldots , J \}$.

Let $\mathbb{N}$ be the set of nonnegative integers, i.e.,
$\mathbb{N} = \{0, 1, 2, \ldots \}$ and let $\Z$ be the set of all
integers, i.e., $\Z = \{\ldots , -2, -1, 0, 1, 2, \ldots \}$.
Without loss of generality, in this paper, we represent a table by
a vector of counts ${\mathbf n}=(n_1,\ldots,n_k)$. Under this
point of view, a contingency table ${\mathbf n}$ can be regarded
as a function ${\mathbf n}: \mathcal{X} \longrightarrow
\mathbb{N}$, but it can also be viewed as a vector ${\mathbf n}
\in \mathbb{N}^k$.

The fiber of an observed table ${\mathbf n}_{\mathrm{obs}}$ with
respect to a function $T:{\mathbb N}^k \longrightarrow {\mathbb
N}^s$ is the set
\begin{equation}
{\mathcal F}_T({\mathbf n}_{\mathrm{obs}}) = \left\{ {\mathbf n} \
| \ {\mathbf n} \in \mathbb{N}^k \ , \ T({\mathbf n}) = T({\mathbf
n}_{\mathrm{obs}})\right\} \, .
\end{equation}
When the dependence on the specific observed table is irrelevant,
we will write simply ${\mathcal F}_T$ instead of ${\mathcal
F}_T({\mathbf n}_{\mathrm{obs}})$.

In mathematical statistics framework, the function $T$ is usually
the minimal sufficient statistic of some statistical model and the
usefulness of enumeration of the fiber ${\mathcal F}_T({\mathbf
n}_{\mathrm{obs}})$ follows from classical theorems such as the
Rao-Blackwell theorem, see e.g. \cite{shao:98}.

When the function $T$ is linear, it can be extended in a natural
way to an homomorphism from ${\mathbb R}^n$ in ${\mathbb R}^s$,
$T$ is represented by an $s \times k$-matrix $A_T$, and its generic
element $A_T(\ell,h)$ is
\begin{equation}
A_T(\ell,h) = T_{\ell}(h) ,
\end{equation}
where $T_\ell$ is the $\ell$-th component of the function $T$. In
terms of the matrix $A_T$, the fiber ${\mathcal F}_T$ can be
easily rewritten in the form:
\begin{equation}
{\mathcal F}_T = \left\{ {\mathbf n} \ | \ {\mathbf n} \in
\mathbb{N}^k \ , \ A_T({\mathbf n}) = A_T({\mathbf
n}_{\mathrm{obs}})\right\} \, .
\end{equation}

To navigate inside the fiber ${\mathcal F}_T$, i.e., to connect
any two tables of the fiber ${\mathcal F}_T$ with a path of
nonnegative tables, algebraic statistics suggests an approach
based on the notion of Markov moves and Markov bases. A Markov
move is any table $\mathbf{m}$ with integer entries that preserves
the linear function $T$, i.e. $T({\mathbf n} \pm {\mathbf m}) =
T({\mathbf n})$ for all ${\mathbf n} \in {\mathcal F}_T$.

A finite set of moves ${\mathcal M}=\{ \mathbf{m}_1, \ldots,
\mathbf{m}_r \}$ is called a {\em Markov basis} if it is possible to
connect any two tables of ${\mathcal F}_T$ with moves in
${\mathcal M}$. More formally, for all ${\mathbf n}_1$ and
${\mathbf n}_2$ in ${\mathcal F}_T$, there exist a sequence of
moves $\{ {\mathbf m}_{i_1}, \ldots , {\mathbf m}_{i_A}\}$ and a
sequence of signs $\{\epsilon_{i_1}, \ldots, \epsilon_{i_A}\}$
such that
\begin{equation} \label{def_mb1}
{\mathbf n}_2 = {\mathbf n}_1 + \sum_{a=1}^A \epsilon_{i_a}
{\mathbf m}_{i_a}
\end{equation}
and
\begin{equation} \label{def_mb2}
{\mathbf n}_1 + \sum_{j=1}^a \epsilon_{i_j} {\mathbf m}_{i_j} \geq
0 \ \ \ \textrm{ for all } \ a =1, \ldots , A \, .
\end{equation}
See \cite{diaconis-sturmfels} for further details on Markov bases.
Given a Markov basis, the Diaconis-Sturmfels algorithm for
sampling from a distribution $\sigma$ on ${\mathcal F}_T$ starts
from a table ${\mathbf n} \in {\mathcal F}_T$ and proceeds at each
step as follows:
\begin{itemize}
\item Choose a move ${\mathbf m} \in {\mathcal M}$ and a sign
$\epsilon= \pm 1$ with probability $1/2$ each independently on
${\mathbf m}$;

\item Generate a random number $u$ from the uniform distribution
${\mathcal U}[0,1]$;

\item If ${\mathbf n}+ \epsilon {\mathbf m} \in {\mathcal F}_T$
and $\min\{ \sigma({\mathbf n}+\epsilon {\mathbf
m})/\sigma({\mathbf n}),1\}>u$, then the Markov chain moves from
the current table ${\mathbf n}$ to ${\mathbf n}+ \epsilon{\mathbf
m}$; otherwise, it stays at ${\mathbf n}$.
\end{itemize}

To actually compute Markov bases, we associate to the problem two
distinct polynomial rings. First, we define $\mathbb{R}[{\mathbf
x}]=\mathbb{R}[x_1, \ldots, x_k]$, i.e., we associate an
indeterminate $x_h$ to any cell of the table; then, we define
$\mathbb{R}[{\mathbf y}] = \mathbb{R}[y_1, \ldots, y_s]$, with an
indeterminate $y_\ell$ for any component of the linear function
$T$. In the following we will use some facts from commutative
algebra, to be found in, e.g., \cite{cox|little|oshea:92}.

The simplest method to compute Markov bases uses the elimination
algorithm:

\begin{itemize}
\item For each column of the matrix $A_T$, define the polynomial
\begin{equation} \label{poly}
f_h = x_h - \prod_{\ell=1}^s y_\ell^{A_T(\ell,h)} \ \mbox{ for } \
h= 1, \ldots, k \, ;
\end{equation}
Then, consider the ideal generated by the polynomials $f_1,
\ldots, f_k$:
\begin{equation} \label{first-id}
{\mathcal I}= \langle f_1, \ldots, f_k \rangle
\end{equation}
in the polynomial ring ${\mathbb R}[{\mathbf x},{\mathbf y}]$;

\item Eliminate the ${\mathbf y}$'s indeterminates, and obtain the
ideal
\begin{equation} \label{toricid}
{\mathcal I}_{A_T}= \mathrm{Elim}({\mathbf y},{\mathcal I})
\end{equation}
in the polynomial ring ${\mathbb R}[{\mathbf x}]$. The ideal
${\mathcal I}_{A_T}$ in Equation \eqref{toricid} is by definition
the toric ideal associated to $A_T$;

\item A Gr\"obner basis of ${\mathcal I}_{A_T}$ is formed by
binomials. Each binomial defines a move of a Markov basis taking
the exponents. Namely, the correspondence between the binomials
and the moves is given by the log-transformation
\begin{equation} \label{logrel}
\log({\mathbf x}^a - {\mathbf x}^b) = a-b \in {\mathbb R}^k \, .
\end{equation}
\end{itemize}

Although faster algorithms have been implemented to compute toric
ideals, the elimination-based algorithm is the simplest one and we
will use this technique in some of the proofs. For details on
computational methods for toric ideals, see
\cite{bigatti|lascala|robbiano:99} and the implementation in {\tt 4ti2}
\citep{4ti2}.

As noted in e.g. \cite{Rapallo07} and
\cite{chen-dinwoodie-dobra-huber2005}, when the entries of table
have an upper bound, the classical notion of Markov basis is not
sufficient to connect all the tables in a fiber. In fact, the
fiber in the bounded case:
\begin{equation}
{\mathcal F}_T^{\mathbf b} = \left\{ {\mathbf n} \ | \ {\mathbf n}
\in \mathbb{N}^k \ , \ T({\mathbf n}) = T({\mathbf
n}_{\mathrm{obs}}) \ , \ {\mathbf n} \leq {\mathbf b} \right\}
\end{equation}
is in general smaller than the unrestricted one.

As shown in Sections \ref{univ-inc} and \ref{univ-part} as well as
\cite{Rapallo07}, the constraint ${\mathbf n} \leq {\mathbf b}$
translates into a linear system by introducing dummy counts
$\overline n_1, \ldots, \overline n_k$ with $n_h+ \overline
n_h=b_h$ for all $h=1, \ldots, k$. Therefore, in the presence of
upper bounds of the cell counts, the Markov basis must be computed
through a Universal Gr\"obner basis of the ideal ${\mathcal
I}_{A_T}$.

The procedure to compute a Universal Gr\"obner basis of the ideal
${\mathcal I}_{A_T}$ is fully described in Chapter $7$ of
\cite{sturmfels:96}. Here we summarize the main steps of the
algorithm. Given the matrix $A_T$, its Lawrence lifting is a
matrix $\Lambda(A_T)$ with dimensions $(s + k) \times (2k)$ and
with block representation
\begin{equation} \label{lawlift}
\Lambda(A_T) = \begin{pmatrix} A_T & 0 \\
                                I_k & I_k
               \end{pmatrix} \, ,
\end{equation}
where $0$ is a null matrix with dimensions $s \times k$ and $I_k$
is the identity matrix with dimension $k \times k$.

The Universal Gr\"obner basis of $A_T$ is then computed with the
algorithm below:
\begin{itemize}
\item Define $k$ new indeterminates $\overline x_1, \ldots,
\overline x_k$;

\item Compute a Gr\"obner basis of the toric ideal ${\mathcal
I}_{\Lambda(A_T)}$ in the polynomial ring ${\mathbb R}[{\mathbf
x},{\mathbf{\overline x}}]$, the toric ideal associated to the
Lawrence lifting $\Lambda(A_T)$ of $A_T$;

\item Substitute $\overline x_h=1$ for all $h= 1, \ldots , k$.
\end{itemize}
The interested reader can find all details and the proof of the
correctness of this algorithm in \cite{sturmfels:96}, Chapter 7.
In terms of Markov bases, we state the following definition.

\begin{definition}
A Markov basis computed through a Universal Gr\"obner basis is a
\emph{Universal Markov basis}.
\end{definition}

Recall that a Universal Gr\"obner basis of the toric ideal
${\mathcal I}_{A_T}$ is formed by binomials, while the
corresponding Universal Markov basis is formed by moves, that is
tables with integer entries. A Gr\"obner basis is a polynomial
object, while a Markov basis is a combinatorial object. As
mentioned above, the connection between Gr\"obner and Markov bases
is given in Equation \eqref{logrel}.


The following section is devoted to the computation of Universal
Markov bases in special settings, such as incomplete tables,
bounds acting on a subset of the full sample space, or strictly
positive bounds.

\section{Universal Markov bases and incomplete tables} \label{univ-inc}

The computation of Universal Markov bases is not easy in practice,
especially for two distinct circumstances:
\begin{itemize}
\item The computation of a Universal Markov basis is based on
twice the number of indeterminates than the standard Markov basis;

\item The number of moves of a Universal Markov basis increases
quickly with the dimension of the contingency table.
\end{itemize}

\begin{example} \label{exind}
Let us consider $I \times J$ contingency tables under
independence model. With fixed marginal totals, and without upper
bounds, a Gr\"obner basis is formed by all $2 \times 2$ minors
(see \cite{diaconis-sturmfels}). This fact can be proved
theoretically and does not need symbolic computations.

In this special case we are also able to characterize the
Universal Gr\"obner basis. Combining Algorithm 7.2 and Corollary
14.12 in \cite{sturmfels:96}, the Universal Gr\"obner basis is
formed by all the binomials:
\begin{equation}\label{completegraph}
x_{i_1j_1}x_{i_2j_2} \ldots x_{i_sj_s} - x_{i_2j_1}x_{i_3j_2}
\ldots x_{i_1j_s} \, ,
\end{equation}
where $(i_1,j_1), (j_1,i_2) , \ldots , (j_s,i_1)$ is a circuit in
the complete bipartite graph with $I$ and $J$ vertices.

This implies that the number of moves needed for the Universal
Markov basis increases much faster with respect to the Markov
basis for the unbounded problem. Just to give the idea of such
increase, we present in the following table the number of moves of
the Gr\"obner bases for square $I \times I$ tables for the first
$I$'s.
\begin{center}
\begin{tabular}{c|cccccc}
 & $2$ & $3$ & $4$ & $5$ & $6$ & $7$ \\ \hline
Standard Markov basis & $1$ & $9$ & $36$ & $100$ & $225$ & $441$
\\ \hline
Universal Markov basis & $1$ & $15$ & $204$ & $3,940$ & $113,865$ & $4,027,161$ \\
\hline
\end{tabular}
\end{center}
\end{example}

To overcome this difficulty it is of major interest to have some
results for the theoretical computation of Universal Markov bases.
The first result in this direction that we present in this section
is related to tables with structural zeros (or incomplete tables).

Let ${\mathcal X}_0 \subset {\mathcal X}$ be the set of structural zeros of the
table, let $T'$ be the function $T$ restricted to ${\mathcal X}' =
{\mathcal X} \setminus {\mathcal X}_0$ and let ${\mathcal
I}'_{A_T}$ be the toric ideal associated with $A_{T'}$

\begin{theorem} \label{Univ_struct}
Let ${\mathbf n}$ be a contingency table and let ${\mathcal
F}_T^{\mathbf b}$ be its bounded fiber under the bound ${\mathbf
n} \leq {\mathbf b}$. Let ${\mathcal X}_0$ be the set of
structural zeros. Then a Universal Gr\"obner basis for the ideal
${\mathcal I}'_{A_T}$ is obtained from the Universal Gr\"obner
basis of ${\mathcal I}_{A_T}$ by removing the binomials involving
indeterminates in ${\mathcal X}_0$.
\end{theorem}
\begin{proof}
Using Theorem 7.1 in \cite{sturmfels:96}, the Universal Gr\"obner
basis has the following two properties: $(a)$ it is unique; $(b)$
it is a Gr\"obner basis with respect to all term orderings on
${\mathbb R}[{\mathbf x}]$.

Without loss of generality, let us suppose that the structural
zeros are the first cells, i.e., ${\mathcal X}_0=\{1, \ldots,
k'\}$. The unique Universal Gr\"obner basis is, from property
$(b)$ above, a basis with respect to the elimination term ordering
for the first $k'$ indeterminates. Then, we apply Theorem 4 in
\cite{Rapallo-2006} and the elimination algorithm.

Following the scheme in Equations \eqref{poly} through \eqref{first-id}
with the matrix $\Lambda(A_T)$, we define the polynomials
\begin{equation*}
f_h = x_h - \overline y_h \prod_{\ell=1}^s y_{\ell}^{A_T(\ell,h)}
\ \ \ \mbox{ for } \ h=1, \ldots, k
\end{equation*}
and
\begin{equation*}
f_{k+h}=\overline x_h-\overline y_h \ \ \ \mbox{ for } \ h=1,
\ldots, k \, .
\end{equation*}
The ideal in Equation \eqref{first-id} becomes
\begin{equation*}
{\mathcal I} = \langle f_1, \ldots, f_k, f_{k+1}, \ldots, f_{2k}
\rangle
\end{equation*}
in the polynomial ring ${\mathbb R}[{\mathbf x},\overline{\mathbf
x}, {\mathbf y},\overline{\mathbf y}]$. Therefore, the toric ideal
${\mathcal I}_{\Lambda(A_T)}$ as in Equation \eqref{toricid} is
\begin{equation}
{\mathcal I}_{\Lambda(A_T)}= \mathrm{Elim}( \{ {\mathbf
y},\overline{\mathbf y} \} , {\mathcal I} ) \, .
\end{equation}
When $x_1 , \ldots , x_{k'}$ are indeterminates associated to
structural zeros, the relevant ideal is
\begin{equation*}
{\mathcal I}' = \mathrm{Elim}(\{x_1 , \ldots, x_{k'}\},{\mathcal
I} )
\end{equation*}
and the Universal Gr\"obner basis of ${\mathcal I}'_{A_T}$ is computed
through
\begin{equation*}
\mathrm{Elim}( \{ {\mathbf y},\overline{\mathbf y} \} , {\mathcal
I}' ) = \mathrm{Elim}( \{ {\mathbf y},\overline{\mathbf y}\} ,
\mathrm{Elim}(\{x_1, \ldots, x_{k'}\}, {\mathcal I})) =
\end{equation*}
\begin{equation*}
= \mathrm{Elim}( \{x_1, \ldots, x_{k'}\}, \mathrm{Elim} ( \{
{\mathbf y},\overline{\mathbf y}\} , {\mathcal I})) =
\mathrm{Elim}(\{x_1, \ldots, x_{k'}\}, {\mathcal
I}_{\Lambda(A_T)})
\end{equation*}
and then substituting $\overline x_h=1$ for all $h$. As the
Universal Gr\"obner basis is in particular a basis with respect to
the elimination term ordering for the indeterminates $x_1, \ldots,
x_{k'}$, this proves that to remove the binomials involving $x_1,
\ldots , x_{k'}$ from ${\mathcal I}_{\Lambda(A_T)}$ is equivalent to
compute the Universal Gr\"obner basis for the incomplete table.
\end{proof}

If one has the Universal Markov basis for the complete
configuration, Theorem \ref{Univ_struct} applies easily. In fact,
using the correspondence between moves and binomials, the theorem
above is clearly equivalent to the following:

\begin{corollary} \label{zero_moves}
Let ${\mathbf n}$ be a contingency table and let ${\mathcal
F}_T^{\mathbf b}$ be its bounded fiber under the bound ${\mathbf
n} \leq {\mathbf b}$. Let ${\mathcal X}_0$ be the set of
structural zeros. Then a Universal Markov basis for ${\mathcal
F}_{T'}^b$ is obtained from a Universal Markov basis for
${\mathcal F}_{T}^b$ by removing the moves involving the cells in
${\mathcal X}_0$.
\end{corollary}

\begin{example} \label{exind2}
Let us consider  $4 \times 4$ contingency tables with fixed
marginal totals, as in Example \ref{exind}. Without structural
zeros, the Universal Markov basis is formed by $204$ binomials:
$36$ moves involving $4$ cells: $96$ moves involving $6$ cells:
and $72$ moves involving $8$ cells.

Suppose that the cell $(1,1)$ is a structural zero. This kind of
table is depicted below, where $0$ means a structural zero, while
the symbol $\bullet$ denotes a non-zero cell.
\begin{equation*}
\begin{pmatrix}
0 & \bullet & \bullet & \bullet \\
\bullet & \bullet & \bullet & \bullet \\
\bullet & \bullet & \bullet & \bullet \\
\bullet & \bullet & \bullet & \bullet
\end{pmatrix}
\end{equation*}
From the complete Universal Markov basis we can remove all moves
involving the structural zero. Applying Corollary
\ref{zero_moves}, we remove: $9$ moves involving $4$ cells: $36$
moves involving $6$ cells: and $36$ moves involving $8$ cells. The
Universal Markov basis in this case has $123$ moves.

Suppose now that the whole main diagonal contains structural
zeros, as in the figure below.
\begin{equation*}
\begin{pmatrix}
0 & \bullet & \bullet & \bullet \\
\bullet & 0 & \bullet & \bullet \\
\bullet & \bullet & 0 & \bullet \\
\bullet & \bullet & \bullet & 0
\end{pmatrix}
\end{equation*}
In this situation we remove: $30$ moves involving $4$ cells: $80$
moves involving $6$ cells: and $66$ moves involving $8$ cells.
Finally, the Universal Markov basis has only $28$ moves.
\end{example}

The last example is a prototype for the quasi-independence models.
Now consider $I \times J$ contingency tables with
structural zeros under quasi-independence model.
\cite{aoki-takemura-2005jscs} computed a unique minimum
Markov basis for $I \times J$ contingency tables with
structural zeros under quasi-independence model.
\begin{definition}[\cite{aoki-takemura-2005jscs}]
Let ${\mathcal X} = \{(i,j) \mid 1 \le i \le I, 1 \le j \le J\}$
be the sample space and let ${\mathcal X}' = {\mathcal X}
\setminus {\mathcal X}_0$ be the set of cells that are not
structural zeros. Also let
\begin{equation*}
F_0(S) = \left\{{\mathbf m} \ \vrule  \ \sum_{j = 1}^J m_{ij} =
\sum_{i = 1}^I m_{ij} = 0, \mbox{ and } m_{ij} = 0 \mbox{ for }
(i, j) \not \in {\mathcal X}_0 \right\} \, .
\end{equation*}
A loop (or loop move) of degree $r$ on ${\mathcal X}'$ is an $I
\times J$ integer array $M_r(i_1, \ldots , i_r ; j_1, . . . , j_r
) \in F_0(S)$, for $1 \leq i_1, \ldots , i_r \leq I$, $1 \leq j_1,
\ldots , j_r \leq J$, where $M_r (i_1, \ldots , i_r ; j_1, \ldots
, j_r )$ has the elements
\begin{equation*}
\begin{array}{c}
m_{i_1j_1} = m_{i_2j_2} = \ldots = m_{i_{r-1}j_{r-1}} = m_{i_r j_r} = 1,\\
m_{i_1j_2} = m_{i_2j_3} = \ldots = m_{i_{r-1}j_r} = m_{i_r j_1} = -1,
\end{array}
\end{equation*}
and all other elements are zero. Also the level indices $i_1, i_2,
\ldots $, and $j_1, j_2, \ldots$ are all distinct, i.e.
\begin{equation*}
i_m \not = i_n \mbox{ and } j_m \not = j_ n \mbox{ for all } m
\not = n \, .
\end{equation*}
\end{definition}
Specifically, a degree $2$ loop $M_2(i_1, i_2; j_1, j_2)$ is
called a basic move.

The support of a loop $M_r(i_1, \ldots , i_r ; j_1, . . . , j_r )$
is the set of its non-zero cells. A loop $M_r(i_1, \ldots , i_r ;
j_1, . . . , j_r )$ is called df $1$ if $R(i_1, \ldots , i_r ;
j_1, \ldots , j_r )$ does not contain support of any loop on $S$
of degree $2,\ldots , r - 1$, where $R(i_1, \ldots , i_r ; j_1,
\ldots , j_r ) = \{(i, j ) | i \in \{i_1, \ldots , i_r \}, j \in
\{j_1, \ldots , j_r \} \}$.

\begin{corollary}[\cite{aoki-takemura-2005jscs}]\label{aokitakemura}
The set of df $1$ loops of degree $2, \ldots , \min\{I, J\}$
constitutes a unique minimal Markov basis for $I \times J$
contingency tables with structural zeros under quasi-independence
model.
\end{corollary}

The examples above show that in many cases the computation of
Universal Markov bases for incomplete tables inherits benefit from
complete tables. In terms of computations, an incomplete table has
less cells than the corresponding complete table and therefore an
incomplete table implies the use of a smaller number of
indeterminates. Nevertheless, in a complete table with symmetric
constraints the Markov bases can be characterized theoretically
(e.g., independence model presented here), and in many cases
the symmetry of the combinatorial problem can lead to substantial
simplifications in the symbolic computation (see in particular
\citep{aoki|takemura:05}). Moreover, following Theorem
\ref{Univ_struct}, in the computation of Universal Markov bases
through elimination we do not introduce new polynomials and, therefore,
we do not increase the degree of the moves, as usual
in the unbounded problems (see \cite{Rapallo-2006}).

\begin{example}
As a different example, where  Markov bases are much simpler,
we present a computation for a $2^{3-1}$ fraction of a factorial
design. The use of Markov bases for fractions are useful for
experiments with Poisson-distributed response variable and the
upper bounds are needed when the response variable is Binomial
(see \cite{aoki|takemura:06}). Here we consider the lattice
$\{-1,1\}^3$ for an experiment with $3$ factors $A$, $B$, and $C$.
The fraction defined by the aliasing equation $AB=1$ consists of
$4$ cells:
\begin{equation}
(-1,-1,-1), \ \ \ (-1,-1,1), \ \ \ (1,1,-1), \ \ \ (1,1,1).
\end{equation}
These $4$ points can be viewed as an incomplete three-way table.
Computing with {\tt CoCoA} \citep{cocoa}, the standard Markov basis
for this incomplete table under the complete independence model
(i.e., with the one-way marginal totals fixed), we obtain only one
move, represented by the binomial:
\begin{equation}
x_{-1-1-1}x_{111}-x_{-1-11}x_{11-1} .
\end{equation}
From this computation we note that:
\begin{itemize}
\item In this example the standard Markov basis has only one
polynomial and therefore it is by definition a Universal Markov
basis;

\item The standard Markov basis for the corresponding complete
table with $8$ cells is formed by $9$ quadratic square-free
binomials, and the corresponding Universal Markov basis for the
bounded problem has $20$ binomials:
\begin{verbatim}
-x[-1,1,-1]x[1,-1,1] + x[-1,-1,-1]x[1,1,1],
-x[-1,1,1]x[1,-1,-1] + x[-1,-1,-1]x[1,1,1],
-x[-1,1,-1]x[1,-1,-1] + x[-1,-1,-1]x[1,1,-1],
x[-1,1,1]x[1,1,-1] - x[-1,1,-1]x[1,1,1],
-x[-1,-1,1]x[-1,1,-1] + x[-1,-1,-1]x[-1,1,1],
-x[-1,1,1]x[1,-1,-1] + x[-1,-1,1]x[1,1,-1],
x[-1,1,1]x[1,-1,-1] - x[-1,1,-1]x[1,-1,1],
x[-1,-1,1]x[1,-1,-1] - x[-1,-1,-1]x[1,-1,1],
-x[1,-1,1]x[1,1,-1] + x[1,-1,-1]x[1,1,1],
-x[-1,1,1]x[1,-1,1] + x[-1,-1,1]x[1,1,1],
-x[-1,1,-1]x[1,-1,1] + x[-1,-1,1]x[1,1,-1],
-x[-1,-1,1]x[1,1,-1] + x[-1,-1,-1]x[1,1,1],
-x[-1,1,1]x[1,-1,1]x[1,1,-1] + x[-1,-1,-1]x[1,1,1]^2,
-x[-1,-1,1]x[-1,1,-1]x[1,-1,-1] + x[-1,-1,-1]^2x[1,1,1],
x[-1,-1,1]x[1,1,-1]^2 - x[-1,1,-1]x[1,-1,-1]x[1,1,1],
x[-1,1,1]x[1,-1,-1]^2 - x[-1,-1,-1]x[1,-1,1]x[1,1,-1],
-x[-1,-1,-1]x[-1,1,1]x[1,-1,1] + x[-1,-1,1]^2x[1,1,-1],
x[-1,1,1]^2x[1,-1,-1] - x[-1,-1,1]x[-1,1,-1]x[1,1,1],
-x[-1,1,-1]^2x[1,-1,1] + x[-1,-1,-1]x[-1,1,1]x[1,1,-1],
-x[-1,1,-1]x[1,-1,1]^2 + x[-1,-1,1]x[1,-1,-1]x[1,1,1]
\end{verbatim}
\end{itemize}
Notice that in a Metropolis-Hastings algorithm one can also make
use of the complete Markov basis and then discard the chosen move
at a given step if it modifies a cell with a structural zero. But
the computations for this example show that the use of such a
strategy leads to a slower convergence of the Markov chain to the
stationary distribution. The use of the Markov basis with the
unique applicable move is essential for a correct use of the
Metropolis-Hastings algorithm.
\end{example}

\section{Markov bases for partially bounded tables} \label{univ-part}

While the problem in the previous section has a positive answer,
in this section we present a problem without a theoretical
solution. Nevertheless, we show how to write the relevant symbolic
computations and we describe explicitly some special examples.

When working with bounded contingency tables, it is a common
situation to have some cell counts bounded and other counts
unbounded. Moreover, some bounds can be treated as unessential. In
this section, we consider two-way contingency tables under
independence model.

It is well known that under the marginal totals each cell count
$n_{ij}$ can not exceed $\min\{{ n}_{i+}, {
n}_{+j}\}$, where ${ n}_{i+}$ is the $i$-th row total and
${n}_{+j}$ is the $j$-th column total. Thus, any bound
exceeding such value can be ignored. Now, we know that:
\begin{itemize}
\item With no upper bounds, we need a Markov basis formed by the
basic moves of the form $\left(\begin{matrix}+1 & -1 \\ -1 & +1
\end{matrix}\right)$ for all $2 \times 2$ minors of the table;

\item With an upper bound for each cell count, we need the
Universal Markov basis formed by all the closed circuits in the
complete bipartite graph with $I$ and $J$ vertices, as discussed
in the previous section.
\end{itemize}

Example \ref{exind} shows that the differences between such two
situations are noticeable in terms of number of moves. We can
conjecture that with some cells bounded and other cells without
bounds we will fall into an intermediate situation, with a
Gr\"obner basis formed by all the degree two by two minors and
some other square-free binomials.

As pointed out in the previous section, the bounds on the cell
counts are represented as linear constraints through the two
identity matrices $I_k$ in the Lawrence lifting $\Lambda(A_T)$,
see Equation \eqref{lawlift}. Thus, for the computation of Markov
bases for partially bounded table, we have to remove from the
block $[I_k, I_k]$ of $\Lambda(A_T)$ the rows corresponding to
cells without upper bound.

To show the behavior of Universal Markov bases with partial
bounds, we present here some numerical examples of Markov bases
computed with {\tt CoCoA}.

\begin{example}\label{eg4}
Consider a $3 \times 3$ contingency table under independence
model. With a bound on all the cells, the Universal Markov basis
has $15$ moves: $9$ moves of the form $\left(\begin{matrix}+1 & -1 \\ -1
& +1 \end{matrix}\right)$ for all $2 \times 2$ minors of the table plus
the $6$ moves of degree $3$ below:
\begin{equation*}
{\mathbf m}_1 = \begin{pmatrix} 0 & -1 & +1 \\
                                -1 & +1 & 0 \\
                                +1 & 0 & -1  \end{pmatrix} \, ,
\end{equation*}
\begin{equation*}
{\mathbf m}_2 = \begin{pmatrix} 0 & -1 & +1 \\
                                +1 & 0 & -1 \\
                                -1 & +1 & 0  \end{pmatrix} \, ,
\end{equation*}
\begin{equation*}
{\mathbf m}_3 = \begin{pmatrix} -1 & 0 & +1 \\
                                +1 & -1 & 0 \\
                                 0 & +1 & -1  \end{pmatrix} \, ,
\end{equation*}
\begin{equation*}
{\mathbf m}_4 = \begin{pmatrix} -1 & 0 & +1 \\
                                 0 & +1 &-1 \\
                                +1 & -1 & 0  \end{pmatrix} \, ,
\end{equation*}
\begin{equation*}
{\mathbf m}_5 = \begin{pmatrix} -1 & +1 & 0 \\
                                 0 & -1 & +1 \\
                                +1 & 0 & -1  \end{pmatrix} \, ,
\end{equation*}
\begin{equation*}
{\mathbf m}_6 = \begin{pmatrix} -1 & +1 & 0 \\
                                +1 & 0 & -1 \\
                                0 & -1 & +1  \end{pmatrix} \, .
\end{equation*}
Now we have computed the Universal Markov basis in three different
situations, with different types of bounds:
\begin{itemize}
\item with a bound only on the cell $(1,1)$, the Universal Markov
basis has $10$ moves: the $9$ basic moves and ${\mathbf m}_2$;

\item with a bound on the three cells on the main diagonal, the
Universal Markov basis has $13$ moves: the $9$ basic moves,
plus ${\mathbf m}_1$, ${\mathbf m}_2$, ${\mathbf m}_4$ and
${\mathbf m}_6$;

\item with a bound on the five block-diagonal cells: $(1,1)$,
$(2,2)$, $(2,3)$, $(3,2)$ and $(3,3)$, the Universal Markov basis
has $12$ moves: the $9$ basic moves, plus ${\mathbf m}_1$,
${\mathbf m}_2$ and ${\mathbf m}_4$;

\item with a bound on all cells but the $(1,1)$, the Universal
Markov basis has $13$ moves: the $9$ basic moves, plus
${\mathbf m}_3$, ${\mathbf m}_4$, ${\mathbf m}_5$ and ${\mathbf
m}_6$.
\end{itemize}
\end{example}

\begin{example}[\cite{aoki-takemura-2005jscs}]
Consider $6 \times 6$ contingency tables of the following form:

\begin{equation*}
\begin{pmatrix}
0 & \bullet & \bullet & 0 & 0 & \bullet \\
\bullet & 0 & \bullet & \bullet & 0 & 0 \\
\bullet & \bullet & 0 & 0 & \bullet & 0 \\
0 & 0 & \bullet & 0 & \bullet & \bullet \\
\bullet & 0 & 0 & \bullet & 0 & \bullet \\
0 & \bullet & 0 & \bullet & \bullet & 0
\end{pmatrix} \, .
\end{equation*}
The reduced Gr\"obner basis with the degree reverse
lexicographical ordering consists of three basic moves, 20 degree
3 loops, 10 degree 4 loops, and 3 degree 5 loops. Note that the
loops of degree 4 and 5 are not df 1. On the other hand, all the
20 loops of degree 3 are df 1. Hence by Corollary
\ref{aokitakemura}, the above three basic moves and 20 degree 3
loops constitute the unique minimal Markov basis.
\end{example}

\section{Markov subbases for bounded and incomplete two-way contingency tables} \label{univ-subbases}

Despite the computational advances presented in the previous
sections, there are applied problems where one may never be able
to compute a Markov basis. Models of no-3-way interaction and
constraint matrices of Lawrence type seem to be arbitrarily
difficult, namely if we vary $I$ and $J$ for $(I, J, K)-$tables,
the degree and support of elements in a minimal Markov bases can
be arbitrarily large \citep{deloera05}. In general, the number of
elements in a minimal Markov basis for a model can be
exponentially many.  Thus, it is important to compute a reduced
number of moves which connect all tables instead of computing a
Markov basis. \cite{Ian2008} discussed that in some cases, such as
logistic regression, positive margins are shown to allow a set of
Markov connecting moves that are much simpler than the full Markov
basis. One such example is shown in \cite{Hara09} where a Markov
basis for a multiple logistic regression is computed by the
Lawrence lifting of this basis. In the case of bivariate logistic
regression, \cite{Hara09} showed  a simple subset of the Markov
basis which connects all fibers with a positive sample size for
each combination of levels of covariates. Such connecting sets
were formalized in \cite{Chen06} with the terminology {\em Markov
subbasis}.

In this section we use a sample space indexed as $\{1,
\ldots, k\}$ instead of $\{1, \ldots, I \} \times \{1, \ldots J
\}$ whenever possible, in order to make the formulae easier to
read.

\begin{definition}[\cite{Chen06}]
A Markov subbasis ${\mathcal M}_{A_T,{\mathbf n}_{\mathrm{obs}}}$
for ${\mathbf n}_{\mathrm{obs}} \in \N^k$ and integer matrix $A_T$
is a finite subset of $\ker(A_T)\cap \Z^k$ such that, for each
pair of vectors ${\mathbf u}, \, {\mathbf v} \in {\mathcal F}_T$,
there is a sequence of vectors ${\mathbf m_i} \in {\mathcal
M}_{A_T,{\mathbf n}_{\mathrm{obs}}}, i = 1,\ldots , l$, such that
\[
{\mathbf u} = {\mathbf v} + \sum_{i = 1}^l{\mathbf m_i},
\]
\[
0 \leq {\mathbf v} + \sum_{i = 1}^j{\mathbf m_i}, \, j = 1, \ldots, l.
\]
The connectivity through nonnegative lattice points only is
required to hold for this specific ${\mathbf n}_{\mathrm{obs}}$.
\end{definition}
Note that ${\mathcal M}_{A_T,{\mathbf n}_{\mathrm{obs}}}$ for
every ${\mathbf n}_{\mathrm{obs}} \in \N^k$ and for a given $A_T$
is a Markov basis ${\mathcal M}_{A_T}$ for $A_T$.

In this section we first study Markov subbases ${\mathcal
M}_{A_T,{\mathbf n}_{\mathrm{obs}}}^b$ for any bounded two-way
contingency tables ${\mathbf n}_{\mathrm{obs}} \in \N^k$ with
positive bounds, i.e., no structural zeros, under independence
model.  Then we study Markov subbases ${\mathcal M}_{A_T,{\mathbf
n}_{\mathrm{obs}}}^b$ for any incomplete $I \times J$ contingency
tables ${\mathbf n}_{\mathrm{obs}} \in \N^k$ with positive
margins, i.e., $A_T({\mathbf n}_{\mathrm{obs}}) > 0$, under
independence model.

To analyze these cases we recall some definitions from commutative
algebra:
\begin{itemize}
\item An ideal ${\mathcal I} \subset {\mathbb R}[\mathbf x]$ is
{\em radical} if
\begin{equation*}
\{f \in {\mathbb R}[\mathbf x] \ | \ f^n \in {\mathcal I} \ \mbox{
for some } n \} = {\mathcal I} \, ;
\end{equation*}

\item Let ${\mathcal I}, \ {\mathcal J}\subset {\mathbb R}[\mathbf
x]$ be ideals.  The quotient ideal $({\mathcal I} : {\mathcal J})$
is defined by:
\[
\left({\mathcal I} : {\mathcal J} \right) = \left\{f \in {\mathbb R}[\mathbf x] \ | \
f \cdot {\mathcal J} \subset {\mathcal I} \right\} \, ;
\]


\item Let $Z = \{z_1, \ldots , z_s\} \subset \R^k$.  A lattice $L$
generated by $Z$ is defined:
\begin{equation*}
L = \Z Z.
\end{equation*}
$M \subset \R^k$ is called a lattice basis of $L$ if each element
in $L$ can be written as a linear integer combination of elements
in $M$. Now a lattice basis for $\ker(A_T)$ has the property that
any two tables can be connected by its vector increments if one is
allowed to swing negative in the connecting path (see Chapter 12
of \cite{sturmfels:96} for definitions and properties of  a
lattice basis).
\end{itemize}
The reader can find in \cite{cox|little|oshea:92} more details on
the definitions above.

\begin{theorem}[\cite{Ian2008}] \label{subbasis2}
Suppose ${\mathcal I}_M$ is a radical ideal, and suppose $M$ is a
lattice basis.  Let $p=x_1 \cdots x_{k}$. For each index $\ell$
with $(A_T)_\ell > 0$, let ${\mathcal I}_\ell = \langle
x_h\rangle_{(A_T)_{\ell,h}>0}$ be the monomial ideal generated by
indeterminates for cells that contribute to margin $\ell$. Let
${\mathcal L}$ be the collection of indices $\ell$ with $(A_T
{\mathbf n}_{\mathrm{obs}})_\ell > 0$. Define
\begin{equation*}
{\mathcal I}_{\mathcal L} = \left({\mathcal I}_M : \prod_{\ell \in
{\mathcal L}} {\mathcal I}_{\ell} \, \right).
\end{equation*}
If
\begin{equation}\label{subbasis3}
\left({\mathcal I}_{\mathcal L}:({\mathcal I}_{\mathcal L}:p) \right)= \langle {
1} \rangle
\end{equation}
then the moves in $M$ connect all the tables in ${\mathcal F}_T$.
\end{theorem}

For computing the following examples we have used the software
{\tt Singular} \citep{Greuel2003}.

\begin{example}[Continue from Example \ref{eg4}]\label{eg1}
Consider again $3 \times 3$ tables with fixed row and column sums,
which are the constraints from fixing sufficient statistics in
independence model, and with all bounded cells. This is equivalent
with $3\times 3 \times 2$ tables with constraints $[A, C]$, $[B,
C]$, $[A, B]$ for factors $A$, $B$, $C$, which would arise for
example in case-control data with two factors $A$ and $B$ at three
levels each.

The constraint matrix that fixes row and column sums in a $3\times
3$ table  gives a toric ideal with a ${3 \choose 2} \times {3
\choose 2}$ element Gr\"obner basis. Each of these moves can be
paired with its signed opposite to get 9 moves of $3\times 3
\times 2$ tables that preserve sufficient statistics. This is
equivalent to 9 moves of the form $\left(\begin{matrix}+1 & -1 \\
-1 & +1 \end{matrix}\right)$ for all $2 \times 2$ minors of the
table for $3\times 3$ tables under independence model (see
Example \ref{eg4}). These elements make an ideal with a Gr\"obner
basis that is square-free in the initial terms, and hence the
ideal is radical (Proposition 5.3 of \cite{Sturmfels2002}). Then
applying Theorem \ref{subbasis2} with nine margins of case-control
counts, i.e., this is equivalent to having the positive
constraints on bounds, namely we have non-zero bounds for all
cells, shows that these 9 moves do connect tables with positive
case-control sums.  The full Markov basis has 15 moves. Therefore,
the Markov subbasis for this table is the standard Markov basis
for a $3 \times 3$ table under independence model.
\end{example}

\begin{example}[\cite{Ian2008}] \label{eg2}
Consider now $4 \times 4$ tables with  fixed row and column
sums as in Example \ref{eg1}, and with all bounded cells. Again,
this is equivalent with $4\times 4 \times 2$ tables with
constraints $[A, C]$, $[B, C]$, $[A, B]$ for factors $A$, $B$ and
$C$, with factors $A$ and $B$ at four levels each.

The constraint matrix that fixes row and column sums in a $4\times
4$ table gives a toric ideal with a ${4 \choose 2} \times {4
\choose 2}$ element Gr\"obner basis. Each of these moves can be
paired with its signed opposite to get 36 moves of $4\times 4
\times 2$ tables that preserve sufficient statistics:
\begin{equation*}
\left(
\begin{array}{cccc}
0&0&0&0\\
+&0&-&0\\
0&0&0&0\\
-&0&+&0\\
\end{array}\right)\, ,
\ \ \
\left(
\begin{array}{cccc}
0&0&0&0\\
-&0&+&0\\
0&0&0&0\\
+&0&-&0\\
\end{array}\right) \, .
\end{equation*}
These elements make an ideal with a Gr\"obner basis that is
square-free in the initial terms, and hence the ideal is radical
(Proposition 5.3 of \cite{Sturmfels2002}). Then applying Theorem
\ref{subbasis2} with sixteen margins of case-control counts, i.e.,
this is equivalent to having the positive conditions on bounds,
namely we have non-zero bounds for all cells, shows that these 36
moves do connect tables with positive case-control sums.  The full
Markov basis has 204 moves. Therefore, the Markov subbasis for
this table is the standard Markov basis for a $4 \times 4$ table
with fixed row and column sums fixed without bounds.
\end{example}

In practice, the algorithm in Theorem \ref{subbasis2} is not
feasible for a large number of cells in a table.

From Examples \ref{eg1} and \ref{eg2} it seems that for bounded
two-way tables with row and column sums fixed we only need a
standard Markov basis for two-way tables with row and column sums
fixed if these bounds are positive.  In fact, by the following
theorem, additional elements in a Universal Markov basis are
needed for incomplete tables, i.e., structural zeros.

\begin{theorem}\label{main}
Consider $I \times J$ tables with row and column sums fixed and
with all cells bounded. If these bounds are positive, then a
Markov subbasis for the tables is the standard Markov basis for
$I \times J$ tables with row and column sums fixed without bounds, i.e.,
the set of basic moves of all $2 \times 2$ minors.
\end{theorem}

In order to prove Theorem \ref{main} we need the following
proposition.
\begin{proposition}\label{prop1}
Let ${\mathcal I}_h = \langle x_h, \ \overline{x}_h \rangle $ for
$h = 1,  \ldots , k=IJ$. Then we have:
\begin{equation*}
\prod_{h = 1}^{k}{\mathcal I}_h = \langle  z_1 \cdots z_{k} \ | \
z_j = x_j \mbox{ or } \overline{x}_j \mbox{ for }j = 1, \ldots , k
\rangle \, .
\end{equation*}
\end{proposition}
\begin{proof}
One can prove this proposition by induction on $k$.  For $k = 2$,
one can verify that using {\tt Singular} \citep{Greuel2003}.
Assume $\prod_{h = 1}^{k}{\mathcal I}_h = \langle  z_1 \cdots
z_{k} \ | \ z_j = x_j \mbox{ or } \overline{x}_j \mbox{ for } j =
1, \ldots , k \rangle$ holds. We want to prove that $\prod_{h =
1}^{k+1} {\mathcal I}_h = \langle z_1 \cdots z_{k+1} \ | \ z_j =
x_j \mbox{ or } \overline{x}_j \mbox{ for } j = 1, \ldots , k+1
\rangle$. We have:
\begin{equation*}
\begin{array}{ccl}
\prod_{h = 1}^{k+1} {\mathcal I}_h & = & \left( \prod_{h = 1}^{k} {\mathcal I}_h \right) \cdot \langle x_{k+1}, \overline{x}_{k+1} \rangle \\
& = & \langle  z_1 \cdots z_{k} \ | \ z_j = x_j \mbox{ or } \overline{x}_j \mbox{ for } j =  1,  \ldots , k \rangle \cdot \langle x_{k+1}, \overline{x}_{k+1} \rangle \\
& = & \langle  z_1 \cdots z_{k+1} \ | \ z_j = x_j \mbox{ or } \overline{x}_j \mbox{ for } j =  1, \ldots , k+1 \rangle \, . \\
\end{array}
\end{equation*}
\end{proof}

Let $M$ be the set of vectors such that
\begin{equation*}
M = \{\pm \left(e_{i_1j_1} + e_{i_2j_2} - e_{i_1j_2} - e_{i_2j_1}\right)\},
\end{equation*}
where $e_{ij} = e_{ijk}$ is defined as an integral array with $1$
at the cell $(i, j, 1)$ and $-1$ at the cell $(i, j, 2)$ and $0$
every other cells. Also let
\begin{equation}\label{IM}
{\mathcal I}_M = \langle
x_{i_1j_1}x_{i_2j_2}\overline{x}_{i_1j_2}\overline{x}_{i_2j_1} -
x_{i_1j_2}x_{i_2j_1} \overline{x}_{i_1j_1}\overline{x}_{i_2j_2} \
| \ i_1 \not = i_2, \, j_1 \not = j_2 \rangle .
\end{equation}

\begin{proof}[Proof of Theorem \ref{main}]
Consider the ideal ${\mathcal I}_M$ in Equation \eqref{IM}. Its
Gr\"obner basis is square-free in the initial terms, and hence the
ideal is radical (Proposition 5.3 of \cite{Sturmfels2002}). Since
${\mathcal I}_M$ in Equation \eqref{IM} is radical, we use Theorem
\ref{subbasis2}. Let ${\mathcal I}_A$ be the toric ideal
associated with the constraint matrix of the tables $I \times J
\times 2$  with constraints $[A, C]$, $[B, C]$, $[A, B]$ for
factors $A$, $B$, and $C$. We want to show
\begin{equation*}
\left({\mathcal I}_M : \prod_{i = 1, \ldots I, \ j= 1 \ldots J}
{\mathcal I}_{ij}\right) = {\mathcal I}_A \, ,
\end{equation*}
where ${\mathcal I}_{ij} = \langle x_{ij}, \overline{x}_{ij}
\rangle $ for $i = 1, \ldots , I, \ j = 1, \ldots , J$. Clearly
$\left({\mathcal I}_M : \prod_{i = 1, \ldots I, \ j= 1 \ldots J}
{\mathcal I}_{ij} \right)\subset {\mathcal I}_A $. Thus we want to
show ${\mathcal I}_A \subset \left({\mathcal I}_M : \prod_{i = 1,
\ldots I, \ j= 1 \ldots J} {\mathcal I}_{ij}\right)$.

By Proposition \ref{prop1},
and Equation (5) on page 193 in
\citep{cox|little|oshea:92}, we only have to show
\begin{equation*}
{\mathcal I}_A \subset \left({\mathcal I}_M: z_{11} \cdots
 {z}_{IJ} \right)
\end{equation*}
where $z_{ij} = x_{ij} \mbox{ or } \overline{x}_{ij}$ for $i = 1, \cdots , I$ and $ j =  1, \ldots , J$.

Let $f \in {\mathcal I}_A$. Then
by the definition of the quotient ideal, we only have to show
\begin{equation*}
(z_{11} \cdots z_{IJ}) \cdot f \in {\mathcal I}_M \, .
\end{equation*}
Assume $I \leq J$ without loss of generality. Also if $I < J$, we
can reduce all moves written in the form of \eqref{completegraph}
to $I \times I \times 2$ tables and other columns are zeros. Thus
we consider $I \times I \times 2$ tables. We will prove this by
induction on $I$. For $I = 3$, one can verify that the statement
holds using {\tt Singular} \citep{Greuel2003}. Assume that the
statement holds for some $I - 1 \geq 3$.  We want to show the
statement holds for $I$. By the inductive assumption we can assume
that $s = I$ in Equation \eqref{completegraph}. Let $f =
x_{i_1j_1}x_{i_2j_2} \cdots
x_{i_Ij_I}\overline{x}_{i_2j_1}\overline{x}_{i_3j_2} \cdots
\overline{x}_{i_1j_I} - x_{i_2j_1}x_{i_3j_2} \cdots
x_{i_1j_I}\overline{x}_{i_1j_1}\overline{x}_{i_2j_2} \cdots
\overline{x}_{i_Ij_I}$. By the symmetry on the row and column
operations on the table $I \times I \times 2$, without loss of
generality we assume $f = x_{11}x_{22} \cdots
x_{II}\overline{x}_{21}\overline{x}_{32} \cdots \overline{x}_{1I}
- x_{21}x_{32} \cdots x_{1I}\overline{x}_{11}\overline{x}_{22}
\cdots \overline{x}_{II}$. This is a binomial representation of a
move on $I \times I \times 2$ tables
\begin{equation*}
\begin{array}{l}
\left(
\begin{array}{cccccc}
1 & 0 & \ldots & 0 &0 & -1\\
-1 & 1 & \ldots & 0  & 0 & 0\\
\vdots & \vdots &\vdots &\vdots &\vdots &\vdots \\
0 & 0 & \ldots & -1 & 1 & 0\\
0 & 0 & \ldots & 0 & -1 & 1\\
\end{array}
\right)
\left(
\begin{array}{cccccc}
-1 & 0 & \ldots & 0 &0 & 1\\
1 & -1 & \ldots & 0  & 0 & 0\\
\vdots & \vdots &\vdots &\vdots &\vdots &\vdots \\
0 & 0 & \ldots & 1 & -1 & 0\\
0 & 0 & \ldots & 0 & 1 & -1\\
\end{array}
\right)  = \\
\left(
\begin{array}{cccccc}
1 & 0 & \ldots & 0 &0 & 0\\
0 & 1 & \ldots & 0  & 0 & 0\\
\vdots & \vdots &\vdots &\vdots &\vdots &\vdots \\
0 & 0 & \ldots & 0 & 1 & 0\\
0 & 0 & \ldots & 0 & 0 & 1\\
\end{array}
\right)
\left(
\begin{array}{cccccc}
0 & 0 & \ldots & 0 &0 & 1\\
1 & 0 & \ldots & 0  & 0 & 0\\
\vdots & \vdots &\vdots &\vdots &\vdots &\vdots \\
0 & 0 & \ldots & 1 & 0 & 0\\
0 & 0 & \ldots & 0 & 1 & 0\\
\end{array}
\right)  -
\left(
\begin{array}{cccccc}
0 & 0 & \ldots & 0 &0 & 1\\
1 & 0 & \ldots & 0  & 0 & 0\\
\vdots & \vdots &\vdots &\vdots &\vdots &\vdots \\
0 & 0 & \ldots & 1 & 0 & 0\\
0 & 0 & \ldots & 0 & 1 & 0\\
\end{array}
\right)
\left(
\begin{array}{cccccc}
1 & 0 & \ldots & 0 &0 & 0\\
0 & 1 & \ldots & 0  & 0 & 0\\
\vdots & \vdots &\vdots &\vdots &\vdots &\vdots \\
0 & 0 & \ldots & 0 & 1 & 0\\
0 & 0 & \ldots & 0 & 0 & 1\\
\end{array}
\right) , \\
\end{array}
\end{equation*}
where the first $I \times I$ table is the first level of the table
and the second table is the second level. We claim that
\begin{equation}\label{eq1}
(z_{11} \cdots z_{II})\cdot f = \sum_{(i, j) = (1,\, 2) \ldots ,
(I-1, \,  I)} {\mathbf x}^{U(i, j)}\overline{\mathbf x}^{V(i,
j)}\left(x_{1i}x_{jj}\overline{x}_{1j}\overline{x}_{ji} -
x_{1j}x_{ji} \overline{x}_{1i}\overline{x}_{jj} \right) \, ,
\end{equation}
where
\begin{equation*}
U(i, j) = \begin{cases}
\left(
\begin{array}{cccc}
2 & 1 & \ldots & 1\\
1 & 2 & \ldots & 1\\
\vdots & \vdots &\vdots &\vdots\\
1 & 1 & \ldots & 2\\
\end{array}
\right) - \left(
\begin{array}{ccccc}
1 & 0 &0 &\ldots & 0\\
0 & 1 & 0 &\ldots  & 0\\
0 & 0 &0 &\ldots  & 0\\
\vdots &\vdots &\vdots &\vdots &\vdots\\
0 & 0 &0 &\ldots & 0\\
\end{array}
\right) - w \mbox{ if } i = 1, j = 2\\
\Sigma_{(i', j') = (1, 2), \ldots, (i - 1, j - 1)} U(i', j')
+ (e_{1, j-1} + e_{j-1, i-1}) - (e_{1, i} + e_{j, j})  \mbox{
else}
\end{cases}
\end{equation*}
and
\begin{equation*}
V(i, j) = \begin{cases} \left(
\begin{array}{ccccc}
0 & 0 & \ldots & 0 & 1\\
1 & 0 & \ldots & 0 & 0\\
0 & 1 & \ldots & 0 & 0\\
\vdots & \vdots & \vdots &\vdots &\vdots\\
0 & 0 & \ldots & 1 & 0\\
\end{array}
\right) - \left(
\begin{array}{cccc}
0 & \ldots & 0 & 1\\
0 & \ldots & 0 & 0\\
\vdots & \vdots &\vdots &\vdots\\
0 & \ldots & 1 & 0\\
\end{array}
\right)  + w \mbox{ if } i = I - 1, j = I\\
\Sigma_{(i', j') = (i + 1, j + 1), \ldots , (I - 1, I)} V(i', j')
+ (e_{1, i+1} + e_{j+1, j+1}) - (e_{1, j} + e_{j, i}) \mbox{ else}
\end{cases} ,
\end{equation*}
where $w \in \{0, 1\}^{I\times J}$ such that 
\[
w_{ij} = \begin{cases}
1 & \mbox{if } z_{ij} = \overline{x}_{ij}\\
0 & \mbox{else.}
\end{cases}
\]
 By the construction of each coefficient,
each monomial in each term cancels out except the monomial with a
negative sign in the first term of the sum and the monomial with a
positive sign in the last term of the sum.  Also simple
calculations show that
\begin{equation*}
u_1 := \left(
\begin{array}{cccc}
0 & \ldots & 1 & 0\\
0 & \ldots & 0 & 0\\
\vdots & \vdots &\vdots &\vdots\\
0 & \ldots & 0 & 1\\
\end{array}
\right) + U(I-1, I) =
\left(
\begin{array}{cccc}
2 & 1 & \ldots & 1\\
1 & 2 & \ldots & 1\\
\vdots & \vdots &\vdots &\vdots\\
1 & 1 & \ldots & 2\\
\end{array}
\right) - w
\end{equation*}
and
\begin{equation*}
v_1 := \left(
\begin{array}{cccc}
0 & \ldots & 0 & 1\\
0 & \ldots & 0 & 0\\
\vdots & \vdots &\vdots &\vdots\\
0 & \ldots & 1 & 0\\
\end{array}
\right) + V(I-1, I)
=
\left(
\begin{array}{ccccc}
0 & 0 & \ldots & 0 & 1\\
1 & 0 & \ldots & 0 & 0\\
0 & 1 & \ldots & 0 & 0\\
\vdots & \vdots & \vdots &\vdots &\vdots\\
0 & 0 & \ldots & 1 & 0\\
\end{array}
\right) + w
\end{equation*}
and
\begin{equation*}
u_2:= \left(
\begin{array}{ccccc}
0 & 1 & \ldots & 0 & 0\\
1 & 0 &\ldots & 0 & 0\\
0 & 0 & \ldots & 0 & 0\\
\vdots & \vdots & \vdots & \vdots &\vdots\\
0 & 0 & \ldots & 0 & 0\\
\end{array}
\right) + U(1, 2) =
\left(
\begin{array}{ccccc}
1 & 1 & \ldots & 1& 2\\
2 & 1 & \ldots & 1& 1\\
1 & 2 & \ldots & 1& 1\\
\vdots & \vdots & \vdots &\vdots &\vdots\\
1 & 1 & \ldots & 2 & 1\\
\end{array}
\right) - w
\end{equation*}
and
\begin{equation*}
v_2 := \left(
\begin{array}{ccccc}
1 & 0 & \ldots & 0 & 0\\
0 & 1 &\ldots & 0 & 0\\
0 & 0 & \ldots & 0 & 0\\
\vdots & \vdots & \vdots & \vdots &\vdots\\
0 & 0 & \ldots & 0 & 0\\
\end{array}
\right) + V(1, 2) =
\left(
\begin{array}{ccccc}
1 & 0 & \ldots & 0 & 0\\
0 & 1 &\ldots & 0 & 0\\
0 & 0 & \ldots & 0 & 0\\
\vdots & \vdots & \vdots & \vdots &\vdots\\
0 & 0 & \ldots & 1 & 0\\
0 & 0 & \ldots & 0 & 1\\
\end{array}
\right) + w.
\end{equation*}
Then we notice that
\begin{equation*}
\begin{array}{ll}
\left[\left(
\begin{array}{cccccc}
1 & 1 & \ldots & 1 &1 & 1\\
1 & 1 & \ldots & 1  & 1 & 1\\
\vdots & \vdots &\vdots &\vdots &\vdots &\vdots \\
1 & 1 & \ldots & 1 & 1 & 1\\
1 & 1 & \ldots & 1 & 1 & 1\\
\end{array}
\right) - w\right]
\left[\left(
\begin{array}{cccccc}
0 & 0 & \ldots & 0 &0 & 0\\
0 & 0 & \ldots & 0  & 0 & 0\\
\vdots & \vdots &\vdots &\vdots &\vdots &\vdots \\
0 & 0 & \ldots & 0 & 0 & 0\\
0 & 0 & \ldots & 0 & 0 & 0\\
\end{array}
\right) + w\right] & + \\
\left(
\begin{array}{cccccc}
1 & 0 & \ldots & 0 &0 & 0\\
0 & 1 & \ldots & 0  & 0 & 0\\
\vdots & \vdots &\vdots &\vdots &\vdots &\vdots \\
0 & 0 & \ldots & 0 & 1 & 0\\
0 & 0 & \ldots & 0 & 0 & 1\\
\end{array}
\right)
\left(
\begin{array}{cccccc}
0 & 0 & \ldots & 0 &0 & 1\\
1 & 0 & \ldots & 0  & 0 & 0\\
\vdots & \vdots &\vdots &\vdots &\vdots &\vdots \\
0 & 0 & \ldots & 1 & 0 & 0\\
0 & 0 & \ldots & 0 & 1 & 0\\
\end{array}
\right)  & = \left( u_1 \right) \left( v_1 \right)\\
\end{array}
\end{equation*}
and
\begin{equation*}
\begin{array}{ll}
\left[\left(
\begin{array}{cccccc}
1 & 1 & \ldots & 1 &1 & 1\\
1 & 1 & \ldots & 1  & 1 & 1\\
\vdots & \vdots &\vdots &\vdots &\vdots &\vdots \\
1 & 1 & \ldots & 1 & 1 & 1\\
1 & 1 & \ldots & 1 & 1 & 1\\
\end{array}
\right) - w\right]
\left[\left(
\begin{array}{cccccc}
0 & 0 & \ldots & 0 &0 & 0\\
0 & 0 & \ldots & 0  & 0 & 0\\
\vdots & \vdots &\vdots &\vdots &\vdots &\vdots \\
0 & 0 & \ldots & 0 & 0 & 0\\
0 & 0 & \ldots & 0 & 0 & 0\\
\end{array}
\right)+ w\right] & + \\
\left(
\begin{array}{cccccc}
0 & 0 & \ldots & 0 &0 & 1\\
1 & 0 & \ldots & 0  & 0 & 0\\
\vdots & \vdots &\vdots &\vdots &\vdots &\vdots \\
0 & 0 & \ldots & 1 & 0 & 0\\
0 & 0 & \ldots & 0 & 1 & 0\\
\end{array}
\right)
\left(
\begin{array}{cccccc}
1 & 0 & \ldots & 0 &0 & 0\\
0 & 1 & \ldots & 0  & 0 & 0\\
\vdots & \vdots &\vdots &\vdots &\vdots &\vdots \\
0 & 0 & \ldots & 0 & 1 & 0\\
0 & 0 & \ldots & 0 & 0 & 1\\
\end{array}
\right) & =  \left( u_2 \right) \left( v_2 \right). \\
\end{array}
\end{equation*}
Thus, ${\mathbf x}^{u_1}\overline{\mathbf x}^{v_1} - {\mathbf
x}^{u_2}\overline{\mathbf x}^{v_2}$ equals to the left hand side
in Equation \eqref{eq1}.
\end{proof}

Now we assume that the given margins are positive for bounded $I
\times J$ tables, i.e., we assume that all row and column sums are
positive.  Without loss of generality, we can assume that all
margins are positive because cell counts in rows and/or columns
with zero marginals are necessary zeros and such rows and/or
columns can be ignored in the conditional analysis.

Let ${\mathcal X} = \{(i,j) \mid 1 \le i \le I, 1 \le j \le J\}$
and let ${\mathcal X}_0$ be a non-trivial subset of ${\mathcal
X}$. Recall that ${\mathcal X}_0$ is the set of structural zeros
of the table. For Examples \ref{eg5} and \ref{eg6}, we used
Theorem \ref{subbasis2}.

\begin{example}\label{eg5}
We consider $3 \times 3$ tables under independence model with all
cells bounded. We assume row and column sums are positive. We have
studied in which ${\mathcal X}_0$ the standard Markov basis for $3
\times 3$ tables, i.e., the set of the 9 moves of the form
$\left(\begin{matrix}+1 & -1 \\ -1 & +1 \end{matrix}\right)$ for
all $2 \times 2$ minors of the table, connects these bounded
tables with positive conditions.  If $|{\mathcal X}_0| = 1$ or
$|{\mathcal X}_0| = 2$ then  Equation in \eqref{subbasis3} holds.
Thus, these 9 moves connect bounded tables.  For $|{\mathcal X}_0|
= 3$, if ${\mathcal X}_0 = \{(1,1), (2, 2), (3, 3)\}$ after an
appropriate interchange of rows and columns, i.e. there are 6
patterns of ${\mathcal X}_0$, then Equation in \eqref{subbasis3}
does not hold. Otherwise for other patterns of ${\mathcal X}_0$,
Equation in \eqref{subbasis3} holds. Thus, 9 moves connect bounded
tables. For $|{\mathcal X}_0| > 3$, if ${\mathcal X}_0$ contains $
\{(1,1), (2, 2), (3, 3)\}$ after appropriate interchange of rows
and columns, then Equation in \eqref{subbasis3} does not hold.
Otherwise for other patterns of $S$,  Equation in
\eqref{subbasis3} holds. Thus, these 9 moves connect bounded
tables.  Even with the positive margin assumption, if ${\mathcal
X}_0 = \{(1,1), (2, 2), (3, 3)\}$, then the basic moves do not
connect incomplete contingency tables, i.e., we need the Universal
Markov basis.
\end{example}

\begin{example}\label{eg6}
We also consider $4 \times 4$ tables under independence model with
all cells bounded.  We assume row and column sums are positive.
After an appropriate interchange of rows and columns, if we have
structural zero constraints on all diagonal cells (i.e., cells
with indices in ${\mathcal X}_0 = \{(i, j): i = j$ for $i = 1,
\ldots , I \}$), then Equation in \eqref{subbasis3} does not hold.
\end{example}

Now we consider $I \times J$ contingency tables with only diagonal
elements being structural zeros under assumption of positive
conditions on row and column sums. \cite{aoki-takemura-2005jscs}
showed the following propositions.
\begin{proposition}\label{prop3}
Suppose we have $I \times J$ tables with fixed row and column sums.
A set of basic moves is a Markov subbasis for $I \times J$
contingency tables, $I, \ J \geq 4$, with structural zeros in only
diagonal elements under the assumption of positive marginals.
\end{proposition}

From Examples \ref{eg5}, \ref{eg6}, and Proposition \ref{prop3},
we have the following open problem.

\begin{problem}\label{prob1}
Suppose we have $I \times J$ tables with fixed row and column sums.
What is the necessary and sufficient condition on ${\mathcal X}_0$ so that
a set of basic moves is a Markov subbasis for $I \times J$
contingency tables with structural zeros in ${\mathcal X}_0$ under the
assumption of positive marginals.
\end{problem}

\section{Discussions}
In this paper we have studied Markov bases and Markov subbases for
bounded contingency tables, showing many ways to compute them.
While Theorem \ref{Univ_struct} applies to incomplete tables,
Theorem \ref{main} considers bounded tables with positive bounds.
In particular, Theorem \ref{main} shows that considering two-way
tables under independence model for bounded tables with
strictly positive bounds, then the set of basic moves, which is
much smaller than the Universal Markov basis, connects the fibers
with given margins. Thus, in practice we do not need to compute
the Universal Markov basis.

In order to prove Problem \ref{prob1} we may be able to apply
Theorem \ref{subbasis2} and mimic the proof for Theorem
\ref{main}. If we can solve Problem \ref{prob1} this would be very
useful in practice because we know exactly when we only need the
set of basic moves of all $2 \times 2$ minors for two-way
incomplete contingency tables.

\section*{Acknowledgements}

The authors thank Ian Dinwoodie (Duke University, Durham, NC) for
computational help and Akimichi Takemura (University of Tokyo,
Japan) for suggesting some of references.  Also we would like to
thank Seth Sullivant (North Carolina State University, NC) for
useful comments to improve this paper. Finally,  we would like to
thank referees for useful comments to improve this paper.

\bibliographystyle{decsci}
\bibliography{refs12}

\end{document}